\documentclass{article}
\usepackage[utf8]{inputenc}
\usepackage[english]{babel}
\usepackage{amsthm,amsmath,amssymb}
\usepackage{geometry}
\newtheorem{theorem}{Theorem}[section]
\newtheorem{corollary}{Corollary}[theorem]
\newtheorem{proposition}{Proposition}[theorem]
\newtheorem{lemma}[theorem]{Lemma}
\newtheorem{mydef}{Definition}
\usepackage{amsmath}
\usepackage{amsfonts}
\usepackage{graphicx}
\usepackage[colorinlistoftodos]{todonotes}
\usepackage[colorlinks=true, allcolors=blue]{hyperref}

\begin{document}
\providecommand{\keywords}[1]{\textbf{Keywords:} #1}
\providecommand{\amsclas}[1]{\textbf{AMS Subject Classifications:} #1}
\title{\textbf{Central Limit Theorems for a Stationary Semicircular Sequence in Free Probability}}
\author{Zhichao Wang\\
Department of Mathematics, Texas A$\&$M University\\
College Station, Tx 7843-3368\\
\href{mailto:wangzc@math.tamu.edu}{wangzc@math.tamu.edu} }
\date{December 10, 2017}
\maketitle              

\begin{abstract}
In this paper, we focus on studying central limit theorems for functionals of some specific stationary random processes. In classical probability theory, it is well-known that for non-linear functionals of stationary Gaussian sequences, we can get a central-limit result via Hermite polynomials and the diagram formula for cumulants. In this paper, the main result is an analogous central limit theorem, in a free probability setting, for real-valued functionals of a stationary semicircular sequence with long-range dependence, namely the correlation function of the underlying time series tends to zero as the lag goes to infinity.  
\end{abstract}
\\
\\
\amsclas{46L54 $\cdot$ 60F05 $\cdot$ 60G10}
\section{Introduction}

The central limit theorem is a crucial result, with a lot of applications, in probability theory and statistics. The traditional central limit theorem is for independent and identically distributed random variables (or freely independent and identically distributed random variables in free probability theory). However, we want to ask what if we do not have independence for random variables. Naturally, the sequence of random variables should be stationary and have the long-range dependence, that is, with a correlation which decays to zero as the lag tends to infinity with a certain rate. Besides, if the random variables are independent, it is trivial that the functionals of the random variables also satisfies the classical central limit theorem. However, if the sequence of random variables is long-range dependent, then this central-limit result become nontrivial anymore. Thus, we mainly studied the central limit theorem on a sequence of nonlinear functionals of stationary and long-range dependent random variables with identical distributions. Many people studied this problem in classical probability. The main results are from L.Graitis and D. Surgailis in \cite{Giraitis1}, and Péter Breuer and Péter Major in \cite{breuer}. The primary result we employed is the theorem 5 in \cite{Giraitis1}. In this paper, Theorem 5 refers in particular to the theorem 5 in \cite{Giraitis1}.

Theorem 5 states a central limit theorem for a nonlinear functional $H(x)$ acted on stationary Gaussian random variables $\{X_{j}\}$ with long-large dependence. In other words, the covariances \[r(t)=Cov(X_{j},X_{j+t})\rightarrow 0 \ (\forall j\in\mathbb{N})\] as $t\rightarrow\infty.$ Then, this theorem solved the asymptotic property in distribution of \[S_{N}=N^{-1/2}\sum_{j=1}^{N}H(X_{j}),\]
when $N\rightarrow\infty$. Generally speaking, the proof of Theorem 5 mainly used the diagram formula in \cite{avram:fox}, the properties of orthogonal polynomials, Hermite polynomials, and the characteristics of partitions on the table in the diagram formula. The diagram formula for cumulants of Hermite polynomials of time series is the key in this proof. We will show this proof precisely again. 

In this paper, we aim at investigating whether a similar central-limit result holds in a free probability setting. Our original results, theorem 4.1 and theorem 4.5 in this paper, is the extension of Theorem 5 into free probability theory. We prove an analogous central-limit result for a stationary and long-range dependent semicircular sequence with a real-valued continuous function of calculus for the semicircular elements. Our idea come from the correspondences between classical probability and free probability. For example, the Gaussian law is changed into the semicircular law in a free probability setting. As an analogue, Hermite polynomials are replaced by Chebyshev polynomials. The notion of Hermite rank is replaced by the Chebyshev rank. Moreover, we can get an analogous diagram formula in a non-commutative probability space, although the partitions must be non-crossing. In fact, it turns out that the approaches used in the Theorem 5 in classical probability can be modified and employed in the free probability setting to prove an analogous central limit theorem. Our results are formulated in terms of the rate of decay of the covariance function, which can also be represented by Chebyshev rank. We will mainly discuss these modifications in section 4. In addition, coincidentally, Ivan Nourdin and Murad S.Taqqu proved a very similar theorem in \cite{ivan}. There are two differences between their result and ours. First, in \cite{ivan}, they only considered the sequence of polynomials of semicircular elements, while we also considered more general case of a sequence of real-valued continuous functions of calculus for semicircular elements, which means the functional has an infinite series expansion in terms of Chebyshev polynomials of second kind. Second, the approaches of this paper are different from those in \cite{ivan}. In this paper, our proofs rely on cumulants and diagram formulas on Chebyshev polynomials. However, in \cite{ivan}, the authors applied a transfer principle established in \cite{nou}.

In the end, we will introduce some results and approaches from Florin Avram and Robert Fox in \cite{avram:fox}, which are quite different from what we discussed before. They used the spectral functions of cumulants to estimate cumulants and employed the matroid to describe the partitions in the diagram formula. These methods are attractive because the conditions in their theorems do not include a lot of approximations, which makes life easier; hence we want to extend these results or methods to the free probability case. However, there is still a long way to go since we need to overcome the difficulties caused by the non-commutative property.
\section{Preliminaries}
In this section, we review briefly, for the readers' convenience, some backgroud knowledge on free probability and some useful methods used in the following sections.
\subsection{Cumulants}
In probability theory, people always use moments to depict distributions of random variables. The importance of moments can be found in the proof of the ordinary central limit theorem. However, cumulants are more useful in our case. The moments and cumulants can be determined by each other and cumulants of a probability distribution provide an alternative to the moments of the distribution.  

Let $\textbf{X}=(X_{1},X_{2},...,X_{n})$ be a random vector. For $\textbf{t}=(t_{1},...,t_{n})\in\mathbb{R}^{n}$, the join moment generating function is \[g_{\textbf{X}}(t_{1},...,t_{n})=\mathbb{E}\left[e^{\textbf{t}\cdot\textbf{X}}\right].\] 
\begin{mydef}
The joint cumulant generating function of random vector $\textbf{X}$ is \[c_{\textbf{X}}(t_{1},...,t_{n}):=\log\left(\mathbb{E}\left[e^{\textbf{t}\cdot\textbf{X}}\right]\right).\]
The joint cumulant of the random variables $(X_{1},X_{2},...,X_{n})$ is
\begin{equation}
\left \langle X_{1},X_{2},...,X_{n} \right \rangle:=\frac{\partial^{n}c_{\textbf{X}}(t_{1},...,t_{n})}{\partial t_{1}\partial t_{2}...\partial t_{n}}{\Big |} _{t_{1}=t_{2}=..=t_{n}=0.}
\end{equation}
\end{mydef}
An important property of joint cumulants is multilinearity:
\begin{equation}
\left \langle X_{1}+Y,X_{2},...,X_{n} \right \rangle=\left \langle X_{1},X_{2},...,X_{n} \right \rangle+\left \langle Y,X_{2},...,X_{n} \right \rangle.
\end{equation}
\begin{mydef}
For a random variable $X$, the k-th cumulant of $X$ is defined by
\begin{equation}
\kappa_{k}(X)=\underset{k}{\underbrace{\left \langle X,X,...,X\right \rangle}}.
\end{equation}
\end{mydef}
Cumulants are also called semi-invariants in some literature. The importance of cumulants comes from the observation that many properties of random variables can be better represented by cumulants than by moments. For instance, the cumulants of a random variable $Z\ \sim \ {\mathcal {N}}(\mu ,\,\sigma ^{2})$ are \[\kappa_{n}(Z)=\left\{\begin{matrix}
 0, & n>2\\ 
\mu, & n=1\\
\sigma^{2},&n=2.
\end{matrix}\right.\]
The cumulants of a random variable can be determined by the moments of the random variable and vice versa. Thus, when proving the central limit theorems, we can also estimate the cumulants of all orders instead of using moments. The advantage of employing cumulants is that, for  $n>2$, we only need to prove that the cumulants $\kappa_{n}$ converge to zero, which makes the estimation easier sometimes. This is the reason why we use the cumulant to prove central limit theorems in this paper rather than the moment. In free probability theory, we can define free cumulants and this concept is very useful in the proof of free central limit theorem with a similar reason discussed above (see \cite{lectures}).
\subsection{The Diagram Formula}
In this paper, we mainly focus on studying non-linear functionals of a time series. The main approach to simplifying the non-linear functionals is orthogonal polynomials, such as Hermite polynomials and Chebyshev polynomials. We want to give a more general definition rather than only introduce some specific orthogonal polynomial sequences. For more details, see \cite{gpecc}.
\begin{mydef}
The polynomial $A_{m}(x)$ $(m=0,1,2...)$ is called m-th Appell (or the m-th Wick power, also denoted as $:x^{(m)}:$ ) if it satisfies one of the equivalent conditions:
\[A'_{n}(x)=nA_{n-1}(x) \ n\in\mathbb{N}\]
or\[\sum_{n=0}^{\infty}z^{n}A_{n}(x)/n!=e^{zx}\left(\int_{\mathbb{R}}e^{zx}d\mu(x)\right)^{-1},\]
for some probability measure $\mu$ on the real line.
\end{mydef}
Let $X_{j}$ be a zero mean stationary sequence with all moments finite and \[Y_{n}=\sum_{j=1}^{n}:X_{j}^{(m)}: \ ,\]
where$:X_{j}^{(m)}:$ is  the m-th Wick power of $X_{j}$ or the m-th Appell polynomial of $X_{j}$. In this paper, we often use specific Appell polynomials, Hermite polynomials $\{H_{m}(x)\}$.  In this paper, we want to show the asymptotic behavior of the cumulants of $Y_{n}$ using the diagram formula from \cite{avram:fox}. The diagram formula is a combinatorial expansion for the cumulants with respect to appell polynomials of random variables. By the multilinearity, the R-th cumulant of $Y_{N}$ is given by
\[\kappa_{R}(Y_{n})=\sum_{j_{1}=1}^{n}...\sum_{j_{R}=1}^{n}\left \langle:X_{j_{1}}^{(m)}:,...,:X_{j_{R}}^{(m)}:\right\rangle.\]
The diagram formula tells us how to compute the cumulant on the right-hand side in terms of cumulants of time series $X_{j}$. In general, we give the diagram formula for $\left \langle:X_{j_{1}}^{(m_{1})}:,...,:X_{j_{R}}^{(m_{R})}:\right\rangle.$ Consider the following table $G_{0}$ with entries of random variables.
\begin{equation}    
\begin{array}{cccc}   
X_{j_{1}} & X_{j_{1}}&X_{j_{1}}&... \\  
X_{j_{2}}& X_{j_{2}}&X_{j_{2}}&...\\
        &... &... &...       \\
X_{j_{R-1}} & X_{j_{R-1}}& X_{j_{R-1}}&...  \\
X_{j_{R}} & X_{j_{R}}& X_{j_{R}}&... \\
  \end{array}            
\end{equation}
In the $k$-th row, there are $m_{k}$ random variables $X_{j_{k}}$, for any $1\leqslant k\leqslant R.$
On the table, consider partitions $P$ satisfy the conditions:\\
\textbf{(A)} No sets $t\in P$ is contained in a single row of the table.\\
\textbf{(B)} For each partition of the rows of the table into two disjoint sets, there is a set $t\in P$ containing an element from each of the two sets. This means the partition is connected.\\
\textbf{(C)} All blocks in the partition are in pairs.

We denote the set of all partitions on $G_{0}$ that satisfies condition (A) as $\mathcal{P}_{1}$ and denote the set of all partitions on $G_{0}$ that satisfies (A) and (B) as $\mathcal{P}_{2}$. Then, the diagram formula states that
\begin{equation}
E\left[:X_{j_{1}}^{(m_{1})}:...:X_{j_{R}}^{(m_{R})}:\right]=\sum_{P\in\mathcal{P}_{1}}\prod_{t\in P}cum(t)
\end{equation}
and
\begin{equation}\label{diagra}
\left \langle:X_{j_{1}}^{(m_{1})}:,...,:X_{j_{R}}^{(m_{R})}:\right\rangle=\sum_{P\in\mathcal{P}_{2}}\prod_{t\in P}cum(t).
\end{equation}
Here, $cum(t)$ denotes the cumulants of the collection of random variables in the subset $t.$ For instance, if the subset $t$ contains two $X_{j_{1}}$, one $X_{j_{2}}$ and two $X_{j_{3}}$, then $cum(t)=\langle X_{j_{1}},X_{j_{1}},X_{j_{2}},X_{j_{3}},X_{j_{3}}\rangle$. Therefore, the cumulant of $Y_{n}$ can be simplified as
\begin{equation}
\kappa_{R}(Y_{n})=\sum_{P}\sum_{j_{1}}^{n}...\sum_{j_{2}}^{n}\prod_{t\in P}cum(t)=:\sum_{P}S_{n}(P).
\end{equation}
This diagram formula is very powerful in our project since it gives us a method to simplify the cumulants of appell polynomials of random variables. Moreover, when the sequence $\{X_{j}\}$ consists of standard Gaussian random variables, the set $\mathcal{P}_{1}$ and $\mathcal{P}_{2}$ must also satisfy condition (C). There are more details about the diagram formula in \cite{gira2}. In this case, we mainly consider the diagram formula for Hermite polynomials. In free probability theory, one can also define an analogous Appell polynomial, called free Appell polynomials. In this case, a specific free Appell polynomials, Chebyshev polynomials of second kind $\{U_{m}(x)\}$, are the usual object people are interested in.
\subsection{The Diagram Formula in Free Probability}
In this section, we show some results in \cite{michael} and \cite{lectures} about free probability. The main result is an analogous diagram formula in non-commutative probability space.
\begin{mydef}
A non-commutative probability space ($\mathcal{A},\phi$) consists of a unital algebra $\mathcal{A}$ over $\mathbb{C}$ and a unital linear functional $\phi:\mathcal{A}\rightarrow\mathbb{C}$ such that $\phi(1_{\mathcal{A}})=1$. The elements in $\mathcal{A}$ are called non-commutative random variables.
\end{mydef}
Above is the definition of plain non-commutative probability spaces without adding any other structures. However, sometimes we want to study the probability spaces with some structures below.
For instance, the trace property for $(\mathcal{A},\phi)$ is that $\phi(ab)=\phi(ba)$ for all $a,b\in\mathcal{A}.$
\begin{mydef}
For the non-commutative probability space ($\mathcal{A},\phi$), suppose $\mathcal{A}$ is a $\ast-$algebra which means $\mathcal{A}$ is also endowed with an antilinear $\ast-$operation $a\mapsto a^{\ast}$ s.t. $(a^{\ast})^{\ast}=a$ and $(ab)^{\ast}=b^{\ast}a^{\ast}$ for all $a,b\in\mathcal{A}$. Meanwhile, if we also have that $\phi(a^{\ast}a)\geq 0$, for any $a\in\mathcal{A}$, then we call $(\mathcal{A},\phi)$ a $\ast-$probability space.
\end{mydef}
More often, we prefer to consider $(\mathcal{A},\phi)$ as the $\ast-$probability space, because, in this case, we can have 2-norm for $\mathcal{A}$ as $\|a\|_{2}=\sqrt{\phi(a^{\ast}a)},$ which means we can take advantage of many analytic techniques to solve problems. The self-adjoint random variables satisfy $a^{*}=a$. Sometimes, we need more structures on the non-commutative probability space. For example, the $\mathcal{A}$ is a $C^{*}-$algebra. After knowing the structures of the non-commutative probability space, we can define moments and free cumulants. In this paper, the free joint cumulant of $X_{1},X_{2},...,X_{k}$ is denoted as $R(X_{1},X_{2},...,X_{k})$; the k-the free cumulant of $X$ is denoted as $R_{k}(X).$ The reader is referred to lecture 8-11 in \cite{lectures} and \cite{michael} for a detailed definition. In order to show some results, in \cite{michael}, of computing the free cumulants of free Appell polynomials, at first, we introduce the following notations about partitions.

Let $G_{1}=\{1,2,...,n_{1},n_{1}+1,...n_{1}+n_{2},...,\sum_{i=1}^{k-1}n_{i}+1,...,N\}$, where $N=\sum_{i=1}^{k}n_{i}.$ To compute the joint cumulants, we need to study all the partitions over the set $G_{1}$. Let $\mathcal{P}(N)$ be the lattice of all the partitions over the set $G_{1}$. The subset $\mathcal{NC}(N)$ is denoted as the lattice of non-crossing partitions. There are three special partitions, denoted as
\begin{equation}
\widehat{0}=\{(1),(2),...,(N)\}
\end{equation}
\begin{equation}
\widehat{1}=\{(1,2,...,N)\}
\end{equation}
\begin{equation}
\pi_{n_{1},n_{2},...,n_{k}}=\{(1,2,...,n_{1}),(n_{1}+1,...n_{1}+n_{2}),...,(\sum_{i=1}^{k-1}n_{i}+1,...,N)\}
\end{equation}
The last partition separates the set into k rows. This is similar with the table $G_{0}$ in section 2.2. Denote each row as $\overrightarrow{u}_{i}=(\sum_{j=1}^{i-1}n_{j}+1,...,\sum_{j=1}^{i}n_{j})$, for $1\leq i\leq k$. The first partition is the largest partition and the second one the smallest one. Consider an equivalence class induced by any partition $\pi$, for $i,j\in G_{1}$, $i\overset{\pi}{\sim}j$ if and only if i and j are in the same subset of $\pi$. We denote the intersection and union in the lattices by $\wedge$ and $\vee.$ 
\begin{mydef}
Let $\pi\in\mathcal{P}(N).$
\begin{enumerate}
\item $\pi$ is non-homogeneous (with respect to $\pi_{n_{1},n_{2},...,n_{k}}$), provided that for any $B\in\pi$, there is not $C\in\pi_{n_{1},n_{2},...,n_{k}}$ such that $B\subset C.$
\item
$\pi$ is inhomogeneous (with respect to $\pi_{n_{1},n_{2},...,n_{k}}$), provided that\\ $\pi_{n_{1},n_{2},...,n_{k}}\wedge\pi=\widehat{0}.$
\item
$\pi$ is connected (with respect to $\pi_{n_{1},n_{2},...,n_{k}}$), provided that $\pi\vee\pi_{n_{1},n_{2},...,n_{k}}=\widehat{1}$
\end{enumerate}
\end{mydef}
In our methods, we need to calculate the joint free cumulants of Chebyshev polynomials of non-commutative random variables. Chebyshev polynomials of the second kind are free Appell polynomials, which are defined in \cite{michael}. They are orthogonal with respect to the semicircular distribution. In fact, these are the only orthogonal free Appell polynomials. The following theorem, which can be regarded as a free diagram formula, is a consequence from theorem 3.16 of \cite{michael}. Here, the table $G_{0}$ is replaced with $G_{1}$.
\begin{theorem}
If free Appell polynomials are $\{A_{n}(x)\}_{n=0}^{\infty}$ and $\{X_{n}\}_{n=0}^{\infty}$ are a sequence of non-commutative random variables, then the joint free cumulant is
\begin{equation}
R(A_{n_{1}}(X_{1}),A_{n_{2}}(X_{2}), ... ,A_{n_{k}}(X_{k}))=\sum_{\pi\in\mathcal{P}_{3}}\prod_{V\in\pi}R(X_{i}:i\in V).
\end{equation}
The set $\mathcal{P}_{3}\subset\mathcal{NC}(N)$ is the collection of all connected, non-crossing and non-homogeneous partitions with respect to $\pi_{n_{1},n_{2},...,n_{k}}$ over the set $G_{1}.$ The set V is any block in the partition $\pi.$ The joint free moments is 
\begin{equation}
M(A_{n_{1}}(X_{1}),A_{n_{2}}(X_{2}),...,A_{n_{k}}(X_{k}))=\sum_{\pi\in\mathcal{P}_{4}}\prod_{V\in\pi}R(X_{i}:i\in V).
\end{equation}
The set $\mathcal{P}_{4}\subset\mathcal{NC}(N)$ is the collection of all non-crossing and non-homogeneous partitions with respect to $\pi_{n_{1},n_{2},...,n_{k}}$ over the set $G_{1}.$ The set V is any block in the partition $\pi.$
\end{theorem}
The diagram formula claims that the computation of free cumulant about Appell polynomials can be changed into the computation of the joint cumulants of the original random variables. If the random variables are special enough, the formula can be shown more briefly. For example, if the the free Appell polynomials are Chebyshev polynomials and $\{X_{i}\}$ is a semicircular family, then partitions in $\mathcal{P}_{3}$ and $\mathcal{P}_{4}$ should be pair partitions, i.e. every block $V$ consists of two points. Then, in this case, all partitions are inhomogeneous with respect to $\pi_{n_{1},n_{2},...,n_{k}}$. In fact, this is what we want in section 4.

\section{Results and Methods in Classical Probability}
In this section, our main goal is to show the theorem and powerful methods presented in \cite{Giraitis1} and \cite{breuer}. They studied the central-limit properties of functionals over stationary Gaussian random variables and their main result from \cite{Giraitis1} is the Theorem 5. We will give a detailed proof of Theorem 5. The proof of Theorem 5 gives us a direct way to prove the theorem in non-commutative case.
\subsection{Notations}
Let $\left\{X_{t}\right\}_{t\in \mathbb{N}}$ be a real stationary mean zero Gaussian sequence with covariance $r(t)=\mathbb{E}[X_{i}X_{i+t}]$ $(\forall i\in \mathbb{N})$ such that $r(0)=1$ and $r(t)\rightarrow 0$ when $t\rightarrow \infty $. Any real-valued function $H(x)\in L^{2}\left (\mathbb{R},\frac{1}{\sqrt{2\pi}}e^{-x^{2}}dx\right)$ has a series expansion of Hermite polynomials, which is $$H\left (x\right)=\sum_{k=0}^{\infty}c_{k}H_{k}\left(x\right)$$
and $\sum c_{k}^{2}k!<\infty.$  We always assume that $c_{0}=0$.
\begin{mydef}
Suppose that $H(x)\in L^{2}\left (\mathbb{R},\frac{1}{\sqrt{2\pi}}e^{-x^{2}}dx\right)$ and $H\left (x\right)=\sum_{k=0}^{\infty}c_{k}H_{k}\left(x\right)$. The Hermite rank of $H(x)$ is defined by
\[k_{*}:=min\{k\in\mathbb{Z}:k\geq 0 \ and \ c_{k}\neq 0\}.\]
\end{mydef}
Let $r_{H}\left(i-j\right)=\mathbb{E}[H\left (X_{i}\right) H\left (X_{j}\right)].$
So , we can find that $r_{H}\left(i-j\right)=r_{H}\left(j-i\right)$, which is true for any integer $i,j$. Set
$S_{N,t}=\sum_{s=1}^{[N\cdot t]}H(X_{s})$, where $0\leq t\in\mathbb{R}$ and $[Nt]$ represents the maximal integer smaller or equal than $Nt$. Then, $\ S_{N}=S_{N,1}.$
\subsection{The Statement of Theorem 5}
Use all assumptions and notations in section 3.1. If the correlation function $r_{H}(n)$ satisfies
\begin{equation}\label{eq13}
\sum_{n}\left | r_{H}(n) \right |< \infty
\end{equation}
 and 
 \begin{equation}
 \sigma^{2}=\sum_{n\in \mathbf{Z}}r_{H}(n)\neq 0,
\end{equation}
then, as $N\rightarrow \infty $, 
\begin{equation}
N^{-\frac{1}{2}}S_{N,t}\overset{d}{\rightarrow} \sigma W(t), 
\end{equation}
where $\left \{ W\left(t\right)  \right \}_{t\geq 0}$ is the standard Wiener process.\\
\subsection{A Detailed Proof of Theorem 5}
\begin{proof}[Proof of Theorem 5]
We only need to prove the statement when $t=1$ because if we have proved the situation $t=1$, then we can use the Donsker's invariance principle to make sure$N^{-\frac{1}{2}}S_{N,t}$ converges in distribution to a standard Brownian motion as $ {\displaystyle n\to \infty .} $ Then, in order to show $$N^{-\frac{1}{2}}S_{N}\overset{d}{\rightarrow} \sigma \mathcal{N}(0,1)$$ we only need to prove that the cumulant generating function of  $N^{-\frac{1}{2}}S_{N}$ converges to the cumulant generating function of the standard Guassian random variable, which is $C_{normal}(u)=\frac{1}{2}u^{2}.$ There are three subproblems we need to solve.
\begin{enumerate}
\item Prove that $E[N^{-\frac{1}{2}}S_{N}]=0$ for any $N\in \mathbb{N}$. This is because of the assumption $c_{0}=0$.
\item Prove that $Var(N^{-\frac{1}{2}}S_{N})\rightarrow \sigma ^{2}$, as $N\rightarrow \infty $.
\item Prove that all the cumulants of order $k\geq 3$ of  $\frac{S_{N}}{\sqrt {N}}$ vanish as $N\rightarrow \infty$.
\end{enumerate}
\leftline{\textbf{Step 1.}}
During the proof, we need to use some properties of Hermite polynomials. First of all, we use a specific diagram formula considering that $\left\{X_{i}\right\}$ are Gaussian random variables with mean zero and$\left\{H_{j}\right\}$ are Hermite polynomials (section 2.2). It states that
\begin{equation}
\mathbb{E}[H_{n_{1}}(X_{1})H_{n_{2}}(X_{2})...H_{n_{k}}(X_{k})]=\sum_{\pi\in \Gamma (n_{1}+...+n_{k})}\prod _{(i,j)\in \pi}\mathbb{E}[X_{i}X_{j}].
\end{equation}
In this formula, $\Gamma (n_{1}+...+n_{k})$ is the set of all inhomogeneous partitions whose blocks are only pairs and every vertex is of degree one. $\Gamma (n_{1}+...+n_{k})$ is a set of partition on set $G_{1}$ in section 2.3 and a subset of $\mathcal{P}(N)$. Applying the formula above, we can easily get
\begin{equation}
\mathbb{E}[H_{k}(X_{0})H_{j}(X_{t})]=\delta(k,j)k!(r(t))^{k}.
\end{equation}
Let $k_{*}$ be the Hermite rank of $H(x)\in L^{2}$. We can compute that 
\begin{equation}
r_{H}(t)=r^{k_{*}}(t)\sum_{n\geq k_{*}}c_{n}^{2}n!r^{n-k_{*}}(t).
\end{equation}
\leftline{\textbf{Step 2.}}
We prove the following statement, which gives us an equivalent property of (\ref{eq13}) in terms of Hermite rank $k_{*}.$
\begin{lemma}
\begin{equation}
\sum_{t}|r_{H}(t)|<\infty \Leftrightarrow \sum_{t} |r(t)|^{k_{*}}<\infty.
\end{equation}
\end{lemma}
\begin{proof}
$\Rightarrow:$ We know that $r(t)\rightarrow 0$, when $t\rightarrow\infty.$
Thus, $\exists T>0$ so that for $\forall |t|>T,$ we have 
\[\sum_{k=k_{*}+1}^{\infty}c_{k}^{2}k!|r(t)|^{k-k_{*}}\leq \frac{1}{2}c_{k_{*}}^{2}k_{*}!.\]
Then, \[\begin{split}
|r_{H}(t)|=&|r(t)|^{k_{*}}\left|\sum_{n\geq k_{*}}c_{n}^{2}n!r^{n-k_{*}}(t)\right|\\
\geq&|r(t)|^{k_{*}}\left[c_{k_{*}}^{2}k_{*}!-\sum_{k=k_{*}+1}^{\infty}c_{k}^{2}k!|r(t)|^{k-k_{*}}\right]\\
\geq&\frac{1}{2}|r(t)|^{k{*}}c_{k_{*}}^{2}k_{*}!,
\end{split} \]
when $|t|>T$. Therefore, $\sum_{t}|r(t)|^{k_{*}}<\infty.$\\
$\Leftarrow:$ By Cauchy-Schwartz inequality, the covariance satisfies 
\[|r(t)|=|Cov(X_{0},X_{t})|=|\mathbb{E}(X_{0}X_{t})|\leq \sqrt{\mathbb{E}(X_{0}^{2})\mathbb{E}(X_{t}^{2})}=1.\]Then, for any $k\geq k_{*}$, $|r(t)|^{k}\leq |r(t)|^{k_{*}}$. Hence, when $k\geq k_{*}$,
\[\sum_{t\in\mathbb{Z}}|r(t)|^{k}\leq \sum_{t\in\mathbb{Z}}|r(t)|^{k_{*}}<\infty.\]
Then, because $\sum_{k\geq k_{*}}c_{k}^{2}k!<\infty,$ we can get the following estimation.
\[\begin{split}
\sum_{t} |r_{H}(t)|&=\sum_{t}\left|\sum_{k\geq k{*}}|c_{k}|^{2}k!r(t)^{k}\right|\\
&\leq\sum_{t\in\mathbb{Z}}\sum_{k\geq k{*}}c_{k}^{2}k!|r(t)|^{k_{*}}\\
&=\left(\sum_{t}|r(t)|^{k_{*}}\right)\cdot\left(\sum_{k\geq k_{*}}c_{k}^{2}k!\right)<\infty.
\end{split}\]
\end{proof}
\leftline{\textbf{Step 3.}}
In this step, we estimate the variance of $S_{N}$ as N goes to infinity.
 \[ \begin{split}
Var(S_{N})&=\mathbb{E}[\sum_{i=1}^{N}H(X_{i})]^{2}\\
&=\sum_{1\leq i,j\leq N}\mathbb{E}[H(X_{i})H(X_{j})]\\
&=\sum_{1\leq i,j\leq N}r_{H}(|i-j|)\\
&=\sum_{t=1}^{N}2(N-t)r_{H}(t)+\sum _{i=1}^{N}r_{H}(0)\\
&=2N\sum_{t=1}^{N}r_{H}(t)-2\sum_{t=1}^{N}t\cdot r_{H}(t)+\sum_{i=1}^{N}C_{0}\\
&=N\sum_{t=-N}^{N}r_{H}(t)-2\sum_{t=1}^{N}t\cdot r_{H}(t)
\end{split}\]\\
In the formula above, we define $C_{0}:=r_{H}(0)=\sum_{k=1}^{\infty}c_{k}^{2}k!< \infty$. We give an estimation of the second summation in the last equality above, considering $C_{1}:=\sum |r_{H}(t)|<\infty$. Then, for $\forall \epsilon>0$, $\exists K\in \mathbf{N}$ such that $\sum _{t>K} |r_{H}(t)|<\epsilon$. Hence, when N is sufficiently large, 
\[ \begin{split}
S_{2}(N):=\frac{1}{N}\left|\sum _{t=1}^{N}t\cdot r_{H}(t)\right|&\leq  \frac{1}{N}\left|\sum _{t\leq K}t\cdot r_{H}(t)\right|+\frac{1}{N}\left|\sum _{N\geq t> K}t\cdot r_{H}(t)\right|\\
&\leq \frac{K}{N}\sum_{t}\left|r_{H}(t)\right|+\frac{N}{N}
\sum_{t>K}\left|r_{H}(t)\right|\\
&\leq \frac{KC_{1}}{N}+\epsilon
\end{split}\]\\
Because $\epsilon$ is arbitrary, $S_{2}(N)$ goes to zero as $N$ goes to infinity. Therefore, if $\sigma ^{2}=\sum r_{H}(t)$,  then we can obtain 
 $\frac{Var(S_{N})}{N}\rightarrow \sigma^{2}.$
Meanwhile,  by step 2 together with what we discussed above, for any $N\in\mathbb{N}$, setting $D:=\sum |r(t)|^{k{*}}$, then
\[\begin{split}
Var(N^{-1/2}\sum_{t=1}^{N}\sum_{k>n}c_{k}H_{k}(X_{t}))&\leq \sum _{k>n}\sum_{t=0}^{N}c_{k}^{2}k!|r^{k}(t)|\cdot2 (N-t)/N\\
&=\sum _{k>n}2c_{k}^{2}k!\sum_{t=0}^{N}\frac{N-t}{N}|r(t)|^{k_{*}}\\
&\leq(\sum_{k>n}2c_{k}^{2}k!)\cdot D\rightarrow 0 \ (n\rightarrow\infty).\end{split}\]
Let $S_{N}^{>n}:=N^{-1/2}\sum_{t=1}^{N}\sum_{k> n}c_{k}^{2}H_{k}(X_{t})$ and $S_{N}^{<n}:=N^{-1/2}\sum_{t=1}^{N}\sum_{k\leq n}c_{k}^{2}H_{k}(X_{t})$.  Then $N^{-1/2}S_{N}=S_{N}^{>n}+S_{N}^{<n}$. Then, for any $t$, $$P(|S_{N}^{>n}|>t)<\frac{Var(S_{N}^{>n})}	{t^{2}}\rightarrow 0,$$ which implies $S_{N}^{>n}\rightarrow 0$ in probability as $n$ goes to infinity. Notice that this convergence does not rely on $N$, because of the uniform convergence of variance of $S_{N}^{>n}$ for $N.$ Then, by virtue of the property of convergence in probability, 
\begin{equation}\label{eq19}
S_{N}^{<n}\rightarrow N^{-1/2}S_{N},(n\rightarrow\infty)
\end{equation} in probability and in distribution, uniformly for any $N\in\mathbb{N}.$

Next, we use the Levy's metric which is $$d(X,Y)=d(F_{X},G_{Y}):=\inf\left\{ \epsilon>0:F(x-\epsilon)-\epsilon\leq G(x)\leq F(x+\epsilon)+\epsilon,\forall x\in \mathbb{R} \right \},$$
where $F_{X}$ and $G_{Y}$ are two cumulative distribution function of random variables $X$ and $Y$. So this metric represents the convergence in distribution for random variables.  If we have proven the theorem for finite series, i.e. for each $n\in \mathbb{N}$, 
\begin{equation}
d(S_{N}^{<n},Y_{n})\rightarrow 0
\end{equation}
as $N\rightarrow \infty$ and $Y_{n}=\sigma_{n}\mathcal{N}(0,1)$, where $\sigma _{n}=\sum_{t}\sum _{k=k_{*}}^{n}c_{k}^{2}k!r^{k}(t)$. So, $\sigma_{n}\rightarrow \sigma$, which means $d(Y_{n}, \sigma\mathcal{N}(0,1))\rightarrow 0$ as n goes to infinity. As what we discussed above, for any $\epsilon>0$, there exists $k\in\mathbb{N}$ such that for any $n\geq k,$ $d(Y_{n}, \sigma \mathcal{N}(0,1))<\frac{\epsilon}{3}$ and by (\ref{eq19}), there also exists $M\in\mathbb{N}$ such that $\forall n\geq M,$ $d(N^{-1/2}S_{N}, S_{N}^{<n})<\frac{\epsilon}{3}$. Let $n=\max\{M,k\}.$ Then, $\exists A(n)\in\mathbb{N}$, such that for any $N>A(n),$ $d(S_{N}^{<n},Y_{n})<\frac{\epsilon}{3}.$ Therefore, for any $N>A(n),$ we have
$$d(N^{-1/2}S_{N},\sigma \mathcal{N}(0,1))\leq d(Y_{n}, \sigma \mathcal{N}(0,1))+d(N^{-1/2}S_{N}, S_{N}^{<n}) +d(S_{N}^{<n},Y_{n})<\epsilon.$$
Then we can get the final result of Theorem 5, i.e. $N^{-1/2}S_{N}\Rightarrow\sigma \mathcal{N}(0,1)$ in distribution.

\leftline{\textbf{Step 4.}}
Therefore, it suffices to prove the theorem when the $H(x)$ only has a finite Hermite series. Assume
\begin{equation}
H(x)=\sum_{k=1}^{m}c_{k}H_{k}(x).
\end{equation}
Let $S_{N,n}=\sum_{t=1}^{N}H_{n}(X_{t})$, where $1\leq n\leq m$. We are supposed to compute all the cumulants $\kappa_{k}(S_{N})=\kappa _{k}(\sum_{n=1}^{m}S_{N,n})$ for every  $k\geq 3$. By the relationship between k-th cumulant and joint cumulant, 
\begin{equation}
\kappa _{k}(\sum_{n=1}^{m}S_{N,n})=\sum \left <S_{N,n_{1}},S_{N,n_{2}},...,S_{N,n_{k}}\right>.
\end{equation} Define 
\begin{equation}
J_{N}= \left<S_{N,n_{1}},S_{N,n_{2}},...,S_{N,n_{k}}\right>
\end{equation}
as the joint cumulant and $1\leq n_{i}\leq m$. Thus, for $k\geq 3$, we want to show that $\kappa_{k}(\frac{S_{N}}{\sqrt{N}})=\frac{\kappa_{k}}{N^{k/2}}\rightarrow 0$ as $N\rightarrow \infty$ , so it suffices to show every \begin{equation}
J_{N}=\left<S_{N,n_{1}},S_{N,n_{2}},...,S_{N,n_{k}}\right>=o(N^{k/2})
\end{equation} for any $n_{i}\geq k_{*}\ (1\leq i\leq k)$.
By the multi-linearity of the joint cumulants, 
\begin{equation}
J_{N}=\sum_{t_{1}}\sum_{t_{2}}...\sum_{t_{k}}\left< H_{n_{1}}(X_{t_{1}}),... ,H_{n_{k}}(X_{t_{k}})  \right>.
\end{equation}
Then, we employ the diagram formula presented in section 2.2 to compute the joint cumulants of Hermite polynomials of Gaussian random variables, i.e.
\begin{equation}\label{eq27}
\left< H_{n_{1}}(X_{t_{1}}),... ,H_{n_{k}}(X_{t_{k}})  \right>=\sum_{\gamma}\prod_{(i,j)}Cov[X_{t_{i}},X_{t_{j}}]=\sum_{\gamma}\prod_{(i,j)}r(t_{i}-t_{j}),
\end{equation}
where $\gamma $ is any connected inhomogeneous partition over the table $G$ below. Equation (\ref{eq27}) is a special case of (\ref{diagra}), according to the normal distribution and Hermite polynomials. Let $J_{N}=\sum_{\gamma}J_{N}(\gamma) $, where
\begin{equation}
J_{N}(\gamma )=\sum _{t_{1},t_{2},...,t_{k}=1}^{N}\prod _{1\leq i<j\leq k}r^{l_{ij}}(t_{i}-t_{j}).
\end{equation}
In this formula, the set $\Gamma $=$\left \{\right.$ $ \gamma$ : connected  graphs  s.t.$\forall \left[ (i,j),(i^{'},j^{'})\right]\in \gamma$ , $i\neq i^{'}$ and  all blocks are pairs$\left \}\right.$. Set  $l_{ij}$ to be the number of edges between the i-th and j-th rows of the table G.\\
\begin{equation}    
G=
\left(                 
  \begin{array}{ccc}   
    (1,1) & ... & (1,n_{1})\\  
  (2,1) & ... & (2,n_{2})\\
        & ...&           \\
  (k,1) & ...  &(k,n_{k})\\
  \end{array}
\right)                 
\end{equation}

The degree of every vertex is one, so it is easy to get that $\sum _{j\neq i}l_{ij}=n_{i},1\leq i\leq k$.Now, we want to separate $J_{N}(\gamma)$ into two parts and estimate them respectively. We write $J_{N}(\gamma)=J_{N}^{'}(\gamma)+J_{N}^{''}(\gamma)$, where
$$J_{N}^{'}(\gamma)=\sum_{t_{1},t_{2},...,t_{k}=1\ |t_{i}-t_{j}|<K, \ if \ l_{ij}>0}^{N} \ \prod_{1\leq i<j\leq k}r^{l_{ij}}(t_{i}-t_{j}).$$
We know that the absolute value of $r(t)$ is smaller or equal than 1. So we only need to estimate the number of indexes $t_{i}$ that satisfies the restriction of $J^{'}_{N}(\gamma )$. By definition, the graph is connected, which means it cannot be separated into two disjoint parts, so we can rearrange the order of rows of the graph such that for $\forall 1\leq i \leq k-1,$ we can find some $i^{'}>i$  s.t.  $l_{ii^{'}}>0$. Then let $s_{i}=t_{i}-t_{i^{'}}$, $1\leq i\leq k-1$. Then, we only require $s_{i}\leq K$, which means we make the $J^{'}_{N}(\gamma)$ larger calculating more terms in the summation. Namely,
\[\begin{split}
\left| J_{N}^{'}(\gamma )\right|& \leq \sum_{t_{1},t_{2},...,t_{k}=1\ |t_{i}-t_{j}|<K, \ if \ l_{ij}>0}^{N} 1\\
&\leq \sum _{|s_{i}|\leq K( 1\leq i\leq k-1 ) \ and \ 1\leq t_{k}\leq N}1\\
&\leq N(2K+1)^{k-1}
\end{split}\]\\
In the second inequality, we only require $\left|t_{i}-t_{i^{'}}\right|\leq K$, which means we have considered some $t_{i}$ and $t_{j}$ into the summation even if $\left|t_{i}-t_{j}\right|\leq K$ and $l_{ij}>0$. Meanwhile we have considered all the $t_{i}$ and $t_{j}$ when $l_{ij}=0$. For the last inequality, we only delete some of $t_{i}$ in each i-th row of the graph
$(1\leq i\leq k-1)$, which satisfies $\left|t_{i}-1\right|>K$, i.e. $\left|t_{i}-t_{i^{'}}\right|>K$ when $t_{i^{'}}=1$.\\
Thus,$\left|\frac{ J_{N}^{'}(\gamma )}{N^{k/2}}\right|\leq \frac{N(2K+1)^{k-1}}{N^{k/2}}\rightarrow 0$, as N goes to infinity for $\forall k\geq 3.$ So we get the estimation of the first part.
\leftline{\textbf{Step 5.}}
For the second part $J_{N}^{''}(\gamma)$, by definition, $\left|J^{''}_{N}(\gamma)\right| \leq \sum _{1\leq i<j\leq k}\left|I_{ij}\right|$, where we define, for instance, 
$$I_{12}:=\sum_{t_{1},t_{2},...,t_{k}=1\ \left |t_{1}-t_{2}\right|>K}^{N}\prod _{1\leq i<j\leq k}r^{l_{ij}}(t_{i}-t_{j})$$
 if$ \ l_{12}>0$ and $I_{12}=0,$ if $l_{12}=0$. WLOG, we only calculate the $I_{12}$, the other $I_{ij}$ have the same estimation as below. We introduce some new notations to simplify the expression of the formula. Let $r_{12}(t,s) =\left| r(t-s)\right|$ if $1\leq t,s\leq N, \ \left| t-s\right|>K$ and is equal to zero otherwise. Let $r_{ij}(t,s)=\left| r(t-s)\right|$ if $1\leq t,s\leq N$ for $(i,j)\neq (1,2)
 $,$1\leq i,j \leq k.$ Then
 $$\left|I_{12}\right|\leq\sum _{t_{1},t_{2},...,t_{k}=1}^{N}\prod _{1\leq i<j\leq k}r^{l_{ij}}_{ij}(t_{i},t_{j})=:R$$
Then, we use Holder's inequality, which is 
\begin{equation}
 |\int \prod _{k=1}^{n}f_{k}|\leq \prod _{k=1}^{n}(\int|f_{k}|^{\beta_{k}})^{1/\beta_{k}}
\end{equation}
and $1/\beta_{1}+...+1/\beta_{k}=1$, to estimate the value above again and again, considering the properties that $l_{ij}=l_{ji}$ and $\sum _{j\neq i}\frac{l_{ij}}{n_{i}}=1$. Then we have the following inequalities
 $$R=\sum _{t_{2},...,t_{k}}(\sum_{t_{1}}\prod_{1\leq i<j\leq k}r^{l_{ij}}_{ij}(t_{i},t_{j}))
 \leq \sum _{t_{2},...,t_{k}}(\sum_{t_{1}}\prod_{j=2}^{k}r^{l_{1j}}_{1j}(t_{1},t_{j}))\prod_{2\leq i<j\leq k}r^{l_{ij}}_{ij}(t_{i},t_{j})
 $$
 $$\leq \sum_{t_{2},...,t_{k}} \prod_{j=2}^{k}(\sum_{t_{1}}(r^{l_{1j}}_{1j}(t_{1},t_{j}))^{\frac{n_{1}}{l_{1j}}})^{\frac{l_{1j}}{n_{1}}}\prod_{2\leq i<j\leq k}r^{l_{ij}}_{ij}(t_{i},t_{j}) $$
$$= \sum_{t_{2},...,t_{k}} \prod_{j=2}^{k}(\sum_{t_{1}}r^{n_{1}}_{1j}(t_{1},t_{j}))^{\frac{l_{1j}}{n_{1}}}\prod_{2\leq i<j\leq k}r^{l_{ij}}_{ij}(t_{i},t_{j}) $$
$$\leq \sum_{t_{3}...t_{k}}(\sum_{t_{2}}(\sum_{t_{1}}r_{12}^{n_{1}}(t_{1},t_{2}))^{\frac{l_{12}}{n_{1}}}\prod _{j=3}^{k}r_{2j}^{l_{2j}}(t_{2},t_{j}))\prod_{j=3}^{k}(\sum _{t_{1}}r_{1j}^{n_{1}}(t_{1},t_{j}))^{\frac{l_{1j}}{n_{1}}} \prod _{3\leq i<j\leq k}r^{l_{ij}}_{ij}(t_{i},t_{j}) $$
$$\leq \sum_{t_{3}...t_{k}}(\sum_{t_{2}}(\sum_{t_{1}}r_{12}^{n_{1}})^{\frac{n_{2}}{n_{1}}})^{\frac{l_{12}}{n_{2}}}\prod _{j=3}^{k}(\sum_{t_{2}}r_{2j}^{n_{2}})^{\frac{l_{2j}}{n_{2}}}\prod _{j=3}^{k}r_{2j}^{l_{2j}}(t_{2},t_{j}))\prod_{j=3}^{k}(\sum _{t_{1}}r_{1j}^{n_{1}}(t_{1},t_{j}))^{\frac{l_{1j}}{n_{1}}} \prod _{3\leq i<j\leq k}r^{l_{ij}}_{ij}(t_{i},t_{j}) $$
$$\leq ...\leq \prod _{1\leq i<j\leq k}R_{ij}.$$\\
 
Here, we define $R_{ij}=(\sum_{t_{i}}(\sum_{t_{j}}r_{ij}^{n_{i}}(t_{i},t_{j}))^{\frac{n_{j}}{n_{i}}})^{\frac{l_{ij}}{n_{j}}}$. By symmetry, we also can get $R\leq \prod R_{ji}$. So, $R\leq min( \prod R_{ji}, \prod R_{ij})$. Let $C=\sum _{t}\left| r^{k_{*}}(t)\right|<+\infty$. Then
 $$R\leq \prod_{(i,j)\neq (1,2)}(\sum_{t_{i}}(\sum_{t_{j}}r_{ij}^{n_{i}})^{\frac{n_{j}}{n_{i}}})^{\frac{l_{ij}}{n_{j}}}\cdot (\sum(\sum r_{12}^{n_{1}})^{\frac{n_{2}}{n_{1}}})^{\frac{l_{12}}{n_{2}}}$$
$$=(\sum_{t}(\sum_{s}|r(|t-s|)|^{n_{i}})^{\frac{n_{j}}{n_{i}}})^{\frac{l_{ij}}{n_{j}}}\cdot(\sum_{t:|t-s|\geq K}(\sum_{s} |r(|t-s|)|^{n_{1}})^{\frac{n_{2}}{n_{1}}})^{\frac{l_{12}}{n_{2}}}
$$ 
$$\leq \prod _{(i,j)\neq(1,2)}(\sum_{t=1}^{N}C^{\frac{n_{j}}{n_{i}}})^{\frac{l_{ij}}{n_{j}}}\cdot W_{12}\leq C^{'}N^{\sum_{2\leq i<j\leq k} \frac{l_{ij}}{n_{j}}}W_{12}$$

In the equalities, we give an estimation of the second term as $$W_{12}:=(\sum_{t:|t-s|\geq K}(\sum_{s} |r(|t-s|)|^{n_{1}})^{\frac{n_{2}}{n_{1}}})^{\frac{l_{12}}{n_{2}}}\leq (N(\sum_{t\geq K}|r(t)|^{n_{1}})^{\frac{n_{2}}{n_{1}}})^{\frac{l_{12}}{n_{2}}}$$
$$\leq N^{\frac{l_{12}}{n_{2}}}(\sum_{t\geq K}|r(t)|^{n_{2}})^{\frac{l_{12}}{n_{1}}}$$
Therefore, $R\leq C^{'}N^{\theta} max\left\{(\sum_{t\geq K}|r(t)|^{n_{2}})^{\frac{l_{12}}{n_{1}}}, (\sum_{t\geq K}|r(t)|^{n_{1}})^{\frac{l_{21}}{n_{2}}} \right\}= C^{'}N^{\theta} \epsilon (K)$, where $\epsilon (K)\rightarrow 0$ as $K\rightarrow \infty$ because $\sum _{t}\left| r^{m}(t)\right|<+\infty .$ The last thing is to estimate the index $$\theta=min\left\{\sum_{1\leq i<j\leq k}\frac{l_{ij}}{n_{j}},\sum_{1\leq i<j\leq k}\frac{l_{ij}}{n_{i}}\right\}.$$

Since $\sum_{j\neq i}l_{ij}=n_{i}$ and $l_{ij}=l_{ji}$, we can figure out that $$\sum_{1\leq i<j\leq k}\frac{l_{ij}}{n_{j}}+\sum_{1\leq i<j\leq k}\frac{l_{ij}}{n_{i}}=\sum_{i=1}^{k}\sum_{j\neq i}\frac{l_{ij}}{n_{i}}=k.$$
Thus, $\theta\leq k/2$, which implies that $R\leq C^{'}N^{k/2}\epsilon (K) $. So $\frac{J_{N}^{''}(\gamma)}{N^{k/2}}\rightarrow 0$ as N goes to infinity for all k larger or equal to 3.
\end{proof}

\section{Analogous Central Limit Theorems in Free Probability}
Now, we want to prove an analogous free central limit theorem as a counterpart of Theorem 5 in non-commutative probability spaces. In the non-commutative probability, the Appell polynomials used in this section are the Chebyshev polynomials of the second kind. The normal distribution is replaced with the semicircular distribution, which is
\begin{equation}\label{eq30}
f(x)=\frac{1}{2\pi}\sqrt{4-x^{2}}\mathcal{X}_{(-2,2)}.
\end{equation}
Any non-commutative random variable with a semicircular distribution is called a semicircular element. Let $\left\{X_{i}\right \}_{i\in\mathbb{Z}_{+}}$ be a sequence of stationary semicircular random variables on a non-commutative probability space $(\mathcal{A},\phi)$. Suppose their mean $\phi(X_{i})=0$ and variance $\phi(X_{i}^{2})=1$. Let the covariance $\rho(k-l):=\phi(X_{k}X_{l})=\phi(X_{k+m}X_{l+m})$ $(\forall m\in\mathbb{N})$. Also, we assume this sequence is long-range dependent, i.e. 
\begin{equation}\label{eq31}
\rho(t)\rightarrow 0, \ |t|\rightarrow\infty.
\end{equation} 
By the property of the semicircular family (see lecture 11 in \cite{lectures}), for any $k>2$,
\begin{equation}\label{cumu}
R(X_{i_{1}},X_{i_{2}},...,X_{i_{k}})=0.
\end{equation}
Let $\left\{U_{n}\right\}$ be the Chebyshev polynomials of second kind. By property (\ref{cumu}) and the free diagram in section 2.3, to compute the joint free cumulants and moments, all the partitions we need to consider are non-crossing and non-homogeneous pair partitions. In other words, the inhomogeneous property is equivalent to the non-homogeneous property in this case. So, the joint moment is 
\begin{equation}\label{eq33}
M(U_{n_{1}}(X_{1}), ... ,U_{n_{k}}(X_{k}))=\sum_{ \pi\in \Gamma^{'}(n_{1}+...+n_{k})}\prod_{(i,j)\in\pi} Cov[X_{i},X_{j}].
\end{equation}
The set $\Gamma^{'}(n_{1}+...+n_{k})$ represents the collection of all non-homogeneous and non-crossing partitions on the diagram, where the blocks are in pairs. The covariance $Cov[X_{i},X_{j}]=R(X_{i},X_{j})$.
Similarly, the joint cumulant is
\begin{equation}
R(U_{n_{1}}(X_{1}), ... ,U_{n_{k}}(X_{k}))=\sum_{\pi\in\Gamma^{''}(n_{1}+...+n_{k})}\prod_{(i,j)\in \pi}Cov[X_{i},X_{j}],
\end{equation}
where the set $\Gamma^{''}(n_{1}+...+n_{k})$ represents the collection of all connected, non-homogeneous and non-crossing partitions on the diagram with the blocks that are in pairs. $\Gamma^{'}(n_{1}+...+n_{k})$ and $\Gamma^{''}(n_{1}+...+n_{k})$ are the subset of $\mathcal{P}_{4}$ and $\mathcal{P}_{3}$ respectively, which are defined in section 2.3.
Here, $Cov[X_{i},X_{j}]=\rho(i-j)$. We aim at extending the Theorem 5 to non-commutative probability spaces. We modify some notations for the non-commutative situation and get similar central-limit results for stationary non-commutative random sequences.

The Chebyshev polynomials of the second kind $\{U_{n}\}$ form a sequence of orthogonal polynomials with the weight $\omega(x)=\frac{1}{2\pi}\sqrt {4-x^{2}}$ on the interval $[-2,2].$ They are defined by the following recursion formulas.
$$U_{0}(x)=1, U_{1}(x)=x, xU_{k}(x)=U_{k-1}+U_{k+1}(x).$$ 
The orthogonality of Chebyshev polynomials is
\begin{eqnarray}
{\displaystyle \frac{1}{2\pi}\int _{-2}^{2}U_{n}(x)U_{m}(x){\sqrt {4-x^{2}}}\,dx={\begin{cases}0&n\neq m,\\{1}&n=m.\end{cases}}}.
\end{eqnarray}
So, $\left \| U_{n} \right \|_{(L^{2},\omega(x))}=1,$ for any nonnegative integer $n$. For any real-valued function $U(x)\in L^{2}\left((-2,2),\frac{\sqrt{4-x^{2}}}{2\pi}\right)$, there exists a unique polynomial expansion 
\begin{equation}\label{eq36}
U(x)=\sum_{n=0}^{\infty}c_{n}U_{n}(x)
\end{equation} and 
\begin{equation}\label{eq37}
\left \| U \right \|^{2}_{(L^{2},\omega(x))}=\sum_{n}c^{2}_{n}<\infty.
\end{equation}
Define the Chebyshev rank $k_{*}$ as the smallest $k\in \mathbb{Z}^{+}$ such that the coefficient $c_{k}\neq 0.$ Like the results in section 3, it turns out that the Chebyshev rank is an important quantity for the following proof of the central-limit result. To make sense, we always assume that the coefficient $c_{0}=0$.
 
Our goal is to study some limiting properties of a sequence of real-valued functionals of semicircular elements with long-range dependence. Thus, we need to introduce the function calculus for a semicircular element. But, at first, we consider the a simple case when the real-valued function $U(x)\in L^{2}\left((-2,2),\frac{\sqrt{4-x^{2}}}{2\pi}\right)$ has a finite Chebyshev polynomial series expansion. Then, $U(x)$ is a polynomial and $U(X)$ is well-defined for any non-commutative random variable $X.$ Because $c_{0}=0$, $\phi(U(X))=0$. Let the covariance function be $\rho_{U}(n)=\phi(U(X_{0})U(X_{n}))$ and $S_{N}=\sum_{i=1}^{N}U(X_{i}).$ Following with the methods used in the proof of Theorem 5, we can get the theorem 4.1.
\begin{theorem}In a non-commutative probability space $(\mathcal{A},\phi),$ $\{X_{i}\}_{i\in\mathbb{Z}^{+}}$ is a sequence of stationary standard semicircular elements with distribution (\ref{eq30}) and long-range dependence satisfying (\ref{eq31}). $U(x)$ has a finite Chebyshev polynomial series expansion and the first coefficient $c_{0}$ is zero.
If the covariance function satisfies
\begin{equation}
\sum_{n\in \mathbb{Z}}\left | \rho_{U}(n) \right |< \infty
\end{equation}
 and 
 \begin{equation}
 \sigma^{2}=\sum_{n\in \mathbb{Z}}\rho_{U}(n)\neq 0.
\end{equation}
Then, as $N\rightarrow \infty $, 
\begin{equation}
N^{-\frac{1}{2}}S_{N}\overset{d}{\rightarrow} \sigma\mathcal{S}(0,1)
\end{equation}
where $\mathcal{S}(0,1)$ is the semicircle distribution with mean zero and variance 1.
\end{theorem}
In the proof of theorem 4.1, we would like to clarify the differences between this proof and Theorem 5. For more details about non-commutative probability space, such as convergence in distribution, see \cite{lectures}. Also, we can generalize theorem 4.1 considering the non-commutative stochastic processes $N^{-\frac{1}{2}}S_{N,t}=N^{-\frac{1}{2}}\sum_{i=1}^{[Nt]}U(X_{i}).$ 
\begin{lemma}
\begin{equation}
\rho_{U_{n}}(t)=\phi[U_{n}(X_{0})U_{n}(X_{t})]=(\rho(t))^{n}
\end{equation}
and 
\begin{equation}
\phi[U_{n}(X_{0})U_{m}(X_{t})]=\delta^{n}_{m}(\rho(t))^{n},\ (t\in\mathbb{N}).
\end{equation}
\end{lemma}
\begin{proof}
We can derive the results from formula (\ref{eq33}) directly. Notice that, in this case, the set $\Gamma^{'}(n+n)$ only has one element since the graph is non-homogeneous and non-crossing. If $n\neq m$, then the set $\Gamma^{'}(n+m)=\emptyset$, since the partition are in pairs.
\end{proof}
\begin{lemma}
With the same conditions in theorem 4.1, we can get
\begin{equation}
\rho_{U}(t)=\sum_{k=k_{*}}^{\infty}|c_{k}|^{2}(\rho(t))^{k},\ (t\in\mathbb{N})
\end{equation}
and
\begin{equation}
\sum_{n}\left | \rho_{U}(n) \right |< \infty \Leftrightarrow \sum_{n}|\rho(n)|^{k_{*}}<\infty.
\end{equation}
\end{lemma}
The proof of lemma 4.3 is similar with what we did when proving Theorem 5 because the function $U(x)$ is a polynomial in fact. Now, we begin to prove theorem 4.1.

\begin{proof}[Proof of theorem 4.1]
By the definition and properties of the convergence in distribution in the free probability theory (see \cite{lectures}), it suffices to give an estimation for free cumulants of all orders of $S_{N}$, such that when $N\rightarrow\infty$, the free cumulants go to the free cumulants of $\mathcal{S}(0,1)$. Obviously, the expectation (i.e.the first order cumulant) of $S_{N}$ is always equal to zero. Consider the $\mathcal{R}-$transform of semicircle law $\mathcal{S}(0,1)$, which is regarded as the analogue of a cumulant generating function. The free cumulants of  $\mathcal{S}(0,1)$ of all orders are equal to zero except for the second cumulant, which is equal to 1. So we want to estimate the variance (i.e. the seconder order cumulant) of $S_{N}$ and the k-th free cumulants $(\forall k>2)$ as N goes to infinity separably.

The variance of $S_{N}$, in our case, is the second moment of $S_{N}$. Like what we did in step 3 of section 3, we have 
 \[ \begin{split}
R_{2}(S_{N})&=\phi(S_{N}S_{N})\\
&=\sum_{1\leq i,j\leq N}\phi(U(X_{i})U(X_{j}))\\
&=\sum_{1\leq i,j\leq N}\rho_{U}(i-j)\\
&=N\sum_{t=-N}^{N}\rho_{U}(t)-2\sum_{t=1}^{N}t\cdot \rho_{U}(t).
\end{split}\]
Thus, when $N\rightarrow\infty$, the second cumulant 
\[R_{2}(N^{-\frac{1}{2}}S_{N})\rightarrow\sum_{t\in\mathbb{Z}}\rho_{U}(t)=\sigma^{2}.\]
Next, we begin to approximate the higher-order cumulants, considering that the function $U(x)$ can be expanded in finite many Chebyshev polynomials.
Suppose 
\[ U(x)=\sum_{n=k_{*}}^{K}c_{n}U_{n}(x).\]
The k-th cumulant of $S_{N}$ is defined as $R_{k}(S_{N})=R(S_{N},...,S_{N})$, where the free cumulant $R:\mathcal{A}^{k}\rightarrow\mathbb{R}$ is a multi-linear functionals. Thus,
\[\begin{split}
R_{k}(S_{N})&=R(S_{N},...,S_{N})\\
&=\sum_{t_{1},t_{2},...,t_{k}=1}^{N} \ \sum_{n_{1},n_{2},...,n_{k}=k_{*}}^{K} c_{n_{1}}c_{n_{2}}...c_{n_{k}}R(U_{n_{1}}(X_{t_{1}}),...,U_{n_{k}}(X_{t_{k}}))\\
\end{split}\]
We only need to compute the following term:
\[\begin{split}
&\sum_{t_{1},t_{2},...,t_{k}=1}^{N}R(U_{n_{1}}(X_{t_{1}}),...,U_{n_{k}}(X_{t_{k}}))\\
=&\sum_{t_{1},t_{2},...,t_{k}=1}^{N}\sum_{\pi\in\Gamma^{''}}\ \prod_{(i,j)\in\pi}\rho(t_{i}-t_{j})\\
=&\sum_{\pi\in\Gamma^{''}(n_{1}+...+n_{k})}\
\sum_{t_{1},t_{2},...,t_{k}=1}^{N}\prod_{1\leq i<j\leq k}\rho^{l_{ij}}(t_{i}-t_{j})
\end{split}\]

The set $\Gamma^{''}(n_{1}+n_{2}+...+n_{k})$ is the collection of all non-homogeneous, connected and non-crossing pair partitions on the diagram $G_{1}$. Since the partitions are all in pairs, we cannot get a partition that has a subset with more than one elements in the same row. Namely, the partitions are all inhomogeneous. Thus, each subset of any partition is a connection between two different rows. $l_{ij}$ is denoted as the number of edges between the i-th row and j-th row of the partition $\pi$. For each partition $\pi\in\Gamma^{''}(n_{1}+n_{2}+...+n_{k})$, we denote 
\begin{equation}
J_{N}(\pi)=\sum _{t_{1},t_{2},...,t_{k}=1}^{N}\prod _{1\leq i<j\leq k}\rho^{l_{ij}}(t_{i}-t_{j}).
\end{equation}
The set $\Gamma^{''}(n_{1}+n_{2}+...+n_{k})$ is a subset of $\Gamma(n_{1}+...+n_{k})$, which is the collection of all inhomogeneous and connected pair partitions. So for any $\pi\in \Gamma^{''}(n_{1}+n_{2}+...+n_{k})$, it is also in $\Gamma(n_{1}+n_{2}+...+n_{k})$. Thus, the estimation for $J_{N}(\gamma)$ in theorem 5 is also true for this theorem. So,
$$J_{N}(\pi)=o(N^{\frac{k}{2}}).$$
Then, as $N$ goes to infinity, the k-th cumulant of $S_{N}$ goes to zero, when k larger than 2.
\end{proof}
Now, we want to study the general case when function $U(x)$ has an infinite Chebyshev polynomial expansion. In this case, we have to define random variable $U(X)$ for the non-commutative random variable $X$ in a specific non-commutative probability space. Consider that $\mathcal{A}$ is a $C^{*}-$algebra and $(\mathcal{A},\phi)$ is a $C^{*}-$probability space.  Then, we can introduce continuous function calculus for elements in $\mathcal{A}$. For a more detailed introduction to $C^{*}-$algebras, we refer to Chapter 4 of \cite{rich}.
\begin{theorem}[Theorem 4.1.3 in \cite{rich}]
If $A$ is a self-adjoint element of a $C^{*}-$algebra $\mathcal{A}$ there is a unique continuous mapping $\Phi:C(sp(A))\rightarrow\mathcal{A}$ and $\Phi(f)=:f(A)$ such that
\begin{enumerate}
\item $\|f(A)\|=\|f\|$;
\item $\bar{f}(A)=[f(A)]^{*}$, where $\bar{f}$ is the conjugate complex function;
\item $f(A)$ is normal.
\end{enumerate}
\end{theorem}
Therefore, for any real-valued function defined on the spectrum of a self-adjoint operator $A$, $f(A)$ is also a self-adjoint element. Moreover, the expectation $\phi(f(A))=\int_{\mathbb{R}} fd\mu_{A}.$ In our case, the semicircular element is self-adjoint with the spectrum on the real line. Let $X$ be a standard semicircular element in a $C^{*}-$probability space $(\mathcal{A},\phi)$ and any $U(x)\in C(sp(X))=C([-2,2]).$ Then, $U(X)$ is well-defined and self-adjoint by theorem 4.4. Also, the continuous function $U(x)$ is in $L^{2}[-2,2],$ so so (\ref{eq36}) and (\ref{eq37}) still hold here in $L^{2}$ space with a weight function, which is the semicircular distribution. Thus, the random variable $U(X)=\sum_{k=1}^{\infty}c_{k}U_{k}(X)$ in $L^{2}$ (see section 2.3), where $U_{k}$ is the Chebyshev polynomial. Because $\phi(U(X))=\int_{[-2,2]} U(x)d\mu_{X}(x),$ the coefficient $c_{0}=0.$

Above, we proved the theorem when the function $U(x)$ can be expanded into a finite series. In order to show the theorem in general case, we now introduce the Stieltjes (Cauchy) transform methods. Given a self-adjoint random variable $X$ with respect to a probability measure $\mu _{X}$,  the Stieltjes transform of $X$ is defined by \[\mathcal{S}_{X}:=\int_{\mathbb{R}}\frac{1}{x-z}d\mu_{X}(x),\]for all $z\in\mathbb{C}_{+}:=\{z=x+iy:y>0\}.$ A sequence of self-adjoint random variables $X_{N}$ converges to a self-adjoint element $X$ in distribution, that is the probability measure $\mu_{X_{N}}\rightarrow\mu_{X}$ converges weakly to a probability measure. By the property of Stieltjes transform, this is equivalent to the convergence of Stieltjes transform, i.e. $\mathcal{S}_{X_{N}}(z)\rightarrow\mathcal{S}_{X}(z)$ pointwise for every $z\in\mathbb{C}_{+}$ (see \cite{ste,ran}).

In the following proof of general case, we will employ some properties of convergence in probability in a $W^{*}-$probability space $(\mathcal{A},\phi)$ over a von Neumann algebra $\mathcal{A}.$ We refer to \cite{oeb1} and \cite{oeb2} for more details. 
\begin{mydef}
Let $a$ and $a_{n}$ be operators in $\bar{\mathcal{A}}$ affiliate with $\mathcal{A}.$ $a_{n}\rightarrow a$ in probability if $|a_{n}-a|\rightarrow 0$ in distribution.
\end{mydef}
From the definition of convergence in probability, it follows immediately that, for self-adjoint elements $a_{n}$ and $a$, $a_{n}-a\rightarrow 0$ in distribution if and only if $a_{n}\rightarrow a$ in probability. Finally,we introduce the fact that convergence in probability implies convergence in distribution in the non-commutative setting. The key of the proof is the Stieltjes transform (see the proof in \cite{oeb1}).
\begin{proposition}[Proposition2.20 in \cite{oeb1}]
Let $\{a_{n}\}$ be a sequence of self-adjoint operators affiliated with a W$^{*}-$probability space $(\mathcal{A},\phi)$. If $a_{n}$ converges to a self-adjoint operator $a$ in probability, then $a_{n}\rightarrow a$ in distribution.
\end{proposition}
Then, combining the methods we summarized above and the proof of theorem 4.1, we can conclude the following theorem when the function $U(x)$ has an infinite series expansion in the Chebyshev polynomials. For simplicity, we continue to use the notations in the theorem 4.1.
\begin{theorem}
In a $W^{*}-$probability space $(\mathcal{A},\phi),$ $\{X_{i}\}_{i\in\mathbb{Z}^{+}}$ is a sequence of stationary standard semicircular elements with distribution (\ref{eq30}) and long-range dependence satisfying (\ref{eq31}). $U(x)$ is any continuous function on $[-2,2].$
If
\begin{equation}
\sum_{n\in \mathbb{Z}}\left | \rho_{U}(n) \right |< \infty
\end{equation}
 and 
 \begin{equation}
 \sigma^{2}=\sum_{n\in \mathbb{Z}}\rho_{U}(n)\neq 0.
\end{equation}
Then, as $N\rightarrow \infty $, 
\begin{equation}
N^{-\frac{1}{2}}S_{N}\overset{d}{\rightarrow} \sigma\mathcal{S}(0,1)
\end{equation}
where $\mathcal{S}(0,1)$ is the semicircle distribution with mean zero and variance 1.
\end{theorem}
\begin{proof}
As we discussed before, we can express continuous functions of calculus for each $X_{i}$ as an infinite series expansion $U(X_{i})=\sum_{k=1}^{\infty}c_{k}U_{k}(X_{i})$ in the $L^{2}$ sense. Let $X_{N}:=N^{-1/2}S_{N}$ and $X_{N}=S_{N}^{>n}+S_{N}^{<n}$, where
\[S_{N}^{>n}:=N^{-1/2}\sum_{i=1}^{N}\sum_{k>n}c_{k}U_{k}(X_{i})\] and \[S_{N}^{<n}:=N^{-1/2}\sum_{i=1}^{N}\sum_{k=1}^{k=n}c_{k}U_{k}(X_{i}).\]
Also, $U(X_{i})$, $X_{N}$, $S_{N}^{<n}$ and $S_{N}^{>n}$ are self-adjoint because the semicircular element is self-adjoint.
In the proof of theorem 4.1, for any $n\in\mathbb{N},$ we have proven that $S_{N}^{<n}\rightarrow Y_{n}:=\sigma_{n}\mathcal{S}(0,1)$ in distribution when $N$ goes to infinity, where $\sigma_{n}^{2}=\sum_{t}\sum_{k=k_{*}}^{n}c_{k}^{2}\rho(t)^{k}.$ And $Y_{n}\rightarrow S:=\sigma\mathcal{S}(0,1)$ in distribution, as $n$ goes to infinity. Recall what we did in section 3 step 3. Similarly, we can get 
\[Var(S_{N}^{>n})=Var(X_{N}-S_{N}^{<n})\rightarrow 0\]
as $n$ goes to infinity, for all $N$ uniformly. Therefore, for each $z\in\mathbb{C}_{+}$, the difference between the 	Stieltjes transform of $S_{N}^{>n}$ and the Stieltjes transform of zero is 
\[\begin{split}
& |\mathcal{S}_{S_{N}^{>n}}(z)-\mathcal{S}_{0}(z)|\\
=& |\int\frac{1}{t-z} d\mu(t)+\frac{1}{z}|\\
\leq & |\int\frac{t}{z(t-z)}d\mu(t)|\\
\leq & \int_{\mathbb{R}}\frac{|t|}{|z|\cdot|t-z|}\\
\leq & \frac{1}{|Im(z)|^{2}}\left(\int|t|^{2}d\mu(t)\right)^{1/2}\\
= & \frac{Var(S_{N}^{>n})}{|Im(z)|^{2}}\rightarrow 0, \ (n\rightarrow\infty).
\end{split}
\]
In other words, $S_{N}^{>n}$ converges to zero in distribution as $n\rightarrow\infty$ for every $N$ uniformly. Then, $X_{N}-S_{N}^{<n}$ converges to zero in distribution as $n\rightarrow\infty$. Therefore, by definition 8 and proposition 4.4.1, $S_{N}^{<n}$ converges to $X_{N}$ in distribution as $n\rightarrow\infty$ for each $N$. Then, the difference of their corresponding Stieltjes transforms converges to zero for each $z\in\mathbb{C}_{+}.$

Given any $\epsilon>0,$ and for each $z\in\mathbb{C}_{+}$ fixed, there exists $K\in\mathbb{N}$ so that $\forall n>K$ and $\forall N\in\mathbb{N}$,
\begin{equation}
\left|\mathcal{S}_{S_{N}^{<n}}(z)-\mathcal{S}_{X_{N}}(z)\right|<\epsilon/3.
\end{equation}
Meanwhile, there exists $M\in \mathbb{N}$ such that $\forall n>M$,
\begin{equation}
\left|\mathcal{S}_{Y_{n}}(z)-\mathcal{S}_{S}(z)\right|<\epsilon/3.
\end{equation}
There also exists $N_{0}\in\mathbb{N}$ such that for all $N>N_{0}$ and $n=\max\{M,K\} $
\begin{equation}
\left|\mathcal{S}_{S_{N}^{<n}}(z)-\mathcal{S}_{Y_{n}}(z)\right|<\epsilon/3.
\end{equation}
Therefore, combining these three inequalities, $\forall N>N_{0}$, we have
\[\left|\mathcal{S}_{X_{N}}(z)-\mathcal{S}_{S}(z)\right|<\epsilon,\]
which means that $X_{N}=N^{-1/2}S_{N}$ converges to semicircular random variable $\sigma \mathcal{S}(0,1)$ in distribution when $N$ goes to infinity.
\end{proof}
\section{Another Method in Classical Probability}
In this section, we give a summary of the results in \cite{florin} and \cite{avram:fox}. The graphical method in these paper also engenders a central limit theorem for stationary sequences, even for more general random variables in classical probability. In general, this method emphasizes the properties of the graph in order to reduce the fussy computation like what we did above. As what we discuss in section 2, in probability theory, our main aim is to find a method of computing the order of magnitude of $S_{n}(P)$ under some circumstances. The author of \cite{avram:fox} mainly employed the spectral functions and graphic method to reduce into an optimal breaking problem. 

To present the main results in \cite{avram:fox}, we need to give some notations and definitions at first. The k-th cumulant spectral function of sequence $X_{j}$ is a function $f^{(k)}(x_{1},x_{2},...,x_{k-1})$ such that 
\[
cum(X_{j_{1}},...,X_{j_{k}})=\int_{[0,1]^{k-1}} f^{(k)}(x_{1},x_{2},...,x_{k-1})e^{2\pi i[x_{1}(j_{1}-j_{k})+...+x_{k-1}(j_{k-1}-j_{k})]}dx_{1}...dx_{k-1}.
\]
The k-th cumulant spectral function only has k-1 variables because \[cum(X_{j_{1}},X_{j_{2}}...,X_{j_{k}})=cum(X_{j_{1}-j_{k}},X_{j_{2}-j_{k}}...,X_{0}).\]
For instance, when $X_{j}=\sum_{r=0}^{j}c_{j-r}\xi_{r}$, where $\xi_{r}$ is an i.i.d. sequence and $\sum (c_{r})^{2}<\infty$, then we can the joint cumulant is $cum(X_{j_{1},...,X_{j_{k}}})=d_{k}\sum_{i}c_{j_{1}-i}...c_{j_{k}-i}$, where $d_{k}$ is the k-th cumulant of $\xi_{j}$. Let $c(x)$ be the Fourier transform of the sequence $c_{j}$. Then, we claim that the k-th cumulant spectral function is 
\[f^{(k)}(x_{1},x_{2},...,x_{k-1})=d_{k}c(x_{1})c(x_{2})...c(x_{k-1})c(-x_{1}-x_{2}-...-x_{k-1}).\]
\begin{proof}
We only need to check the formula 
\[\int f^{(k)}(x_{1},x_{2},...,x_{k-1})e^{2\pi i[x_{1}(j_{1}-j_{k})+...+x_{k-1}(j_{k-1}-j_{k})]}dx_{1}...dx_{k-1}=d_{k}\sum_{i}c_{j_{1}-i}...c_{j_{k}-i}.\]
If we consider the properties of Fourier transform, then when fixing $x_{2},...,x_{k-1}$
\[\begin{split}
&\int_{0}^{1}c(x_{1})c(-x_{1}-x_{2}-...-x_{k-1})e^{2\pi ix_{1}(j_{1}-j_{k})}dx_{1}\\
=&\mathcal{F}(c(x_{1}))\ast \mathcal{F}(c(-x_{1}-x_{2}-...-x_{k-1}))\\
=&\sum_{i_{1}}c_{j_{1}-j_{k}-i_{1}}c_{-i_{1}}e^{2\pi i(-x_{2}...-x_{k-1})i_{1}}.
\end{split}\]
Then, 
\[\begin{split}
&\int f^{(k)}(x_{1},x_{2},...,x_{k-1})e^{2\pi i[x_{1}(j_{1}-j_{k})+...+x_{k-1}(j_{k-1}-j_{k})]}dx_{1}...dx_{k-1}\\
=&d_{k}\sum_{i_{1}}c_{j_{1}-j_{k}-i_{1}}\int c(x_{2})...c(x_{k-1})e^{2\pi i(x_{2}(j_{2}-j_{k}-i_{1})+...+x_{k-1}(j_{k-1}-j_{k}-i_{1}))}dx_{2}...dx_{k-1}\\
=&d_{k}\sum_{i_{1}}c_{j_{1}-j_{k}-i{1}}c_{j_{2}-j_{k}-i{1}}...c_{j_{k-1}-j_{k}-i{1}}\\
=&d_{k}\sum_{i}c_{j_{1}-i}...c_{j_{k}-i}.
\end{split}
\]
\end{proof}
Another notation is the greatest cross-norm on the tensor product space $L^{(k)}_{p}$ of $L_{p}$ with k times. If for each k-th cumulant spectral function $f^{(k)}(x_{1},...,x_{k-1})$, $k\geq 2$, there exists a function $g^{(k)}(x_{1},...,x_{k})$ and constant $p_{k}\in [0,\infty]$ such that\\
\textbf{(C)} $f^{(k)}(x_{1},...,x_{k-1})=g^{(k)}(x_{1},...,x_{k-1},-x_{1}-...-x_{k-1})$ and\\
\textbf{(D)} $|\|g^{(k)}\||_{p_{k}}<\infty.$\\
Here $|\|g^{(k)}\||_{p_{k}}$ is the greatest cross norm which is defined
\[|\|g^{(k)}\||_{p_{k}}=inf\sum_{j=1}^{N}\|g_{j,1}\|_{p_{k}}...\|g_{j,k}\|_{p_{k}},\]
where the infimum is taken over all decompositions of $g^{(k)}=\sum_{j=1}^{N}g_{j,1}...g_{j,k}$ ($g_{j,i}(x_{i})\in L^{p}$).
\begin{lemma}
If $f(x_{1},...,x_{k-1})=g(x_{1},...,x_{k-1},-x_{1}-...,x_{k-1}),$ then $\|f\|_{p} \leq |\|g\||_{kp}$.
\end{lemma}
\begin{proof}
For any decomposition of function g, $g(x_{1},...,x_{k})=\sum_{i=1}^{N}a_{1i}(x_{1})...a_{ki}(x_{k}),$ we have
\[
\begin{split}
\|f\|_{p}&=\|\sum_{i=1}^{N}a_{1i}(x_{1})...a_{ki}(-x_{1}-...-x_{k-1})\|_{p}\\
&\leq\sum_{i=1}^{N}\|a_{1i}(x_{1})...a_{ki}(-x_{1}-...-x_{k-1})\|_{p}\\
&\leq\sum_{i}^{N}\|a_{1i}\|_{kp}...\|a_{ki}\|_{kp}.
\end{split}\]
The last inequality is a Holder inequality. We can take the infimum over all decompositions.
\end{proof}
The last new definition we want to introduce is the optimal breaking problem. For each partition $P$ that we introduced in the section 2, we can construct a new graph $G$ with  two types of vertices:
\begin{enumerate}
\item 
The element from the set $R$ of vertices is called row vertex since each vertex represents a row in the table.
\item
The element from the set $T$ of vertices is called subset vertex since each vertex represents a subset in the partition $P.$
\end{enumerate}
The edges of graph $G$ only connect the row vertices and subset vertices, showing that that subsets of partition intersect with some of rows. Namely, each element of the table is represented by an edge connecting row and subset vertices in the graph $G$. On each edge, we associate a cost $z_{k}$, where k is the cardinality of the partition subset containing that element of the table and $z_{k}=1-(p_{k})^{-1}.$ $p_{k}$ is defined above. Let E be the edge set of this graph. We consider a breaking profit for the subset of edges $A\subset E$ as
\begin{equation}
\alpha(A)=C(G\setminus A)-\sum_{e\in A}z_{e},
\end{equation}
where $z_{e}$ is denoted as the cost of edge $e$ and $C(G\setminus A)$ is the number of components left in graph $G$ after the edges in $A$ have been removed. The optimal breaking problem is to find 
\begin{equation}
\alpha_{G}=\underset{A\subset E}{max} \ \alpha(A).
\end{equation}
The paper \cite{avram:fox} showed that the order of magnitude of $S_{n}(P)$, defined in section 2, is $\alpha_{G}.$
\begin{theorem}
Suppose that the cumulant spectral functions of $X_{j}$ satisfy conditions \textbf{(C)(D)} and partition $P$ satisfies conditions \textbf{(A)(B)}. Let $S_{n}(P)$ be the corresponding term in the the expansion of the R-th cumulant of $Y_{n}.$ Then,
\begin{enumerate}
\item 
$|S_{n}(P)|\leq Cn^{\alpha_{G}}$, where C is a constant depending only on the norms $|\|g^{(k)}\||_{p_{k}}$.
\item
If $\alpha_{G}>1$, then $S_{n}(P)=o(n^{\alpha_{G}}).$
\item
If $\alpha_{G}=1$, then $lim_{n\rightarrow\infty}S_{n}(P)/n=I_{G}$, where $I_{G}$ is an integral defined as follows.
\end{enumerate}
Let $t=1,2,..,T$ denote the subsets in the partition $P$ and $n_{t}$ denote the cardinality of subset t. Let the matrix $M^{\ast}$ be an integer representation of the cutset matroid $\mathcal{E}^{\ast}(G)$ of the graph G. Then 
\[I_{G}=\int\prod_{t=1}^{T}f^{(n_{t})}(x_{t,1},...,x_{t,n_{t}})dy_{1}...dy_{N}\]
where the vectors $x=(x_{1,1},...,x_{1,n_{1}},x_{2,1},...,x_{T,1},...,x_{T,n_{T}}) $ and $y=(y_{1},...,y_{N})$ are related by $x=yM^{\ast}$, with $N$ being the number of rows in $M^{\ast}$.
\end{theorem} 
Let $\mathcal{G}_{R}$ be the family of graphs arising from partitions involved in the expansion if the R-th cumulant if $Y_{n}.$
\begin{corollary}
Suppose  that the cumulant spectral functions of $X_{j}$ satisfy conditions \textbf{(C)(D)} and $\alpha_{G}\leq R/2$ for every $G\in \mathcal{G}_{R}$ when $R\geq 2$. Then, $n^{-\frac{1}{2}}Y_{n}$ converges in law to the normal distribution with mean 0 and variance $\sigma^{2}=\sum_{G\in\mathcal{G}_{2}}I_{G}$.
\end{corollary}
This corollary shows in which circumstance we can get a central limit theorem. However, in the result, we cannot say precisely what $\alpha_{G}$ really is. In the following theorem, we can say actually the value of $\alpha_{G}$ relies on both the graph G and the spectral cumulant function. 
\begin{theorem}
If \textbf{(C)} and \textbf{(D)} for the spectral cumulant functions hold, and $z_{k}$ (or $p_{k}$) satisfies
\begin{eqnarray}
1-\frac{1}{p_{k}}=z_{k}\geq {\begin{cases}\frac{k}{2m}&if \ k(k-1)>2m,k\leq m+1,\\\frac{k}{2m(k-m)}&if \ m+1\leq k<2m,\\
\frac{1}{k}+\frac{1}{2m}&otherwise,\end{cases}}
\end{eqnarray}
then $n^{-\frac{1}{2}}Y_{n}$ tends in distribution to the $\mathcal{N}(0,\sigma^{2}).$
\end{theorem}
Using this theorem and corollary, we want to regain what we showed in section 3. In that case, the sequence $X_{j}$ are Gaussian random variables and the Appell polynomials are Hermite polynomials. Then all partitions in this situation satisfy \textbf{(A)(B)}. Moreover, we know all partitions are in pairs. If we want to employ theorem 5.3, then we must assume that
\[1-p_{2}^{-1}=z_{2}\geq \frac{1}{2}+\frac{1}{2m},\]
i.e. $p_{2}\geq \frac{2m}{m-1},$ where m is corresponding with the m-th Hermite polynomials.
Suppose $H(x)\in L^{2}\left (\mathbf{R},e^{-x^{2}}dx \right )$ and the Hermite rank of $H(x)$ is one. Then, we can easily figure out that $H\in L^{\infty}$, hence $p_{2}=\infty$ and we can directly employ theorem 5.3 to get a central-limit result for $\{H(X_{n})\}$. However, the question is, for the rank larger than one, we cannot use these theorems directly. In other words, these theorems in \cite{avram:fox} are weaker than Theorem 5 in section 3 in some sense. Therefore, to apply these methods, we should add more conditions to get the same result. 

In addition, we are motivated by the fact that theorem 5.2 is more concise and general than what we considered in section 3. Hence, we want to extend these theorems in free probability theory. However, we have no idea how to extend these methods to free probability space till now since the spectral method here does not work well in non-commutative case. 

\section*{Acknowledgements}
This paper is a summary of my capstoe project, supervised by Professor Michael Anshelevich. Due to a cooperation math program between Beihang University and Texas A$\&$M University, I was able to gain an invaluable experience over my senior year, studying at Texas A$\&$M University and doing the research with Professor Michael Anshelevich. I thank my advisor, Professor Michael Anshelevich, for his patient guidance. I am also grateful to Beihang University for supporting me to study abroad. 
%


\begin{thebibliography}{5}
%
\bibitem{michael}
Michael Anshelevich: Appell polynomials and their relatives. Int Math Res Notices 2004; 3469-3531 (2009).
\bibitem {avram:fox}
Florin Avram, Robert Fox: Central Limit Theorems for Sums of Wick Products of Stationary Sequences. Transactions of the American Mathematical Society, vol. 330, no. 2, pp. 651–663 (1992).
\bibitem{florin}
Florin Avram: Generalized Szegő theorems and asymptotics of cumulants by graphical methods. Transactions of the American Mathematical Society, Volume 330, Pages 637-649 (1992).
\bibitem{Giraitis1}
 L.Giraitis, D.Surgailis: CLT and other limit theorems for functionals of Gaussian processes. (1985) 70: 191. doi: 10.1007/BF02451428.
\bibitem{gira2}
 L. Giraitis, D. Surgailis: Multivariate Appell polynomials and the central limit theorem, Dependence in Probability and Statistics (Oberwolfach, 1985), Progr. Probab. Statist., vol. 11,
Birkhauser Boston,Massachusetts, pp. 21–71 (1986).
\bibitem{breuer}
Péter Breuer, Péter Major: Central limit theorems for non-linear functionals of Gaussian fields. Journal of Multivariate Analysis, Volume 13, Issue 3, Pages 425-441, ISSN 0047-259X (1983).
\bibitem{lectures}
Alexandru Nica, Roland Speicher, Lectures on the combinatorics of free probability. Cambridge University Press, 2006.
\bibitem{ivan}
I. Nourdin, M.S. Taqqu: Central and non-central limit theorems in a free probability setting. Journal of Theoretical Probability, 27(1), pp.220-248 (2014).
\bibitem{oeb1}
O.E. Barndorff-Nielsen, S. Thorbjørnsen: Self-decomposability and Lévy processes in free probability. Bernoulli, 8(3), pp.323-366 (2002).
\bibitem{oeb2}
O.E. Barndorff-Nielsen, S. Thorbjørnsen: The Lévy-Itô decomposition in free probability. Probability theory and related fields, 131(2), pp.197-228 (2005).
\bibitem{rich}
R. V. Kadison, J. R.Ringrose: Fundamentals of the Theory of Operator Algebras. Vol. I: Elementary Theory.  Academic Press. XV, 398 pp., ISBN 0-12-393301-3 (1983).
\bibitem{ste}
J.S. Geronimo, T.P. Hill: Necessary and sufficient condition that the limit of Stieltjes transforms is a Stieltjes transform. Journal of Approximation Theory, 121(1), pp.54-60 (2003).
\bibitem{nou}
I. Nourdin, G.Peccati, R. Speicher: Multi-dimensional semicircular limits on the free Wigner chaos. In Seminar on Stochastic Analysis, Random Fields and Applications VII (pp. 211-221). Birkhäuser, Basel (2013).
\bibitem{gpecc}
G. Peccati, M. S. Taqqu, Wiener Chaos: Moments, Cumulants and Diagrams, Springer-Verlag, Italia, 2011.
\bibitem{ran}
Michael Anshelevich, Gregory Berkolaiko. Topics in Random matrices Preliminary version \url{http://www.math.tamu.edu/~berko/teaching/spring2017/RMT/Random-matrices-notes-short.pdf} (updated April 27, 2017)
\end{thebibliography}
\end{document}